\documentclass[12pt]{article}
\usepackage{amssymb, url, amsthm, amsmath}
\usepackage{amscd}

\newtheorem{theorem}{\noindent Theorem}

\newtheorem{definition}{\noindent Definition}

\date{}
\title{Dynamics of metrics in measure spaces
and their asymptotic invariants}

\author{A.~M.~Vershik}

\begin{document}
\maketitle
\rightline{On the 50th anniversary of the Kolmogorov entropy}

\begin{abstract}
We discuss the Kolmogorov's entropy and Sinai's definition of it; and then define a deformation of
the entropy, called {\it scaling entropy}; this is also a metric
invariant of the measure preserving actions of the group, which is more powerful than the ordinary entropy.
To define it, we involve the notion of the $\varepsilon$-entropy of a metric in a measure space,
also suggested by A.~N.~Kolmogorov slightly earlier. We suggest to replace
the techniques of measurable partitions, conventional in entropy theory, by that of
iterations of metrics or semi-metrics. This leads us to the key idea of
this paper which as we hope is the answer on the old question: what is the natural context in which one should
consider the entropy of measure-preserving actions of groups? the same question about
its generalizations---scaling entropy, and more general problems of ergodic theory.
Namely, we propose a certain research program, called {\it asymptotic dynamics of
metrics in a measure space}, in which, for instance, the generalized entropy
is understood as {\it the asymptotic Hausdorff dimension of a sequence of metric spaces
associated with dynamical system.} As may be supposed, the metric isomorphism problem
for dynamical systems as a whole also gets a new geometric interpretation.
\footnote{This article is completed mathematical addendum to the paper of the author
"Information, Entropy, Dynamics" in the volume
"Mathematical events of the XX century:look from St.Petersburg". MCCME. 2009.(In Russian),\cite{E}
in which a historical survey of the discovery of the mathematical entropy by Shannon and Kolmogorov
and its influence on mathematics was given.}
\footnote{The author acknowledges support
from RFBR (under grant 08-01-00379) and NSh-2460.2008.1 and thanks
Alexander von Humboldt Foundation (Germany) and Max-Planck-Institute (Bonn)
for providing excellent research facilities.}
\end{abstract}

\section{Introduction}

More than fifty years ago a A.N.Kolmogorov after serious investigation of the notion of Shannon
information theory introduced a new metric invariant in ergodic theory {\it entropy of automorphism of the Lebesgue space}. That event had drastically changed the theory of dynamical systems in all its aspects -  smooth dynamics,
as well as topological and measure-theoretical. After the paper by V.Rokhlin, Ya.Sinai, M.Pinsker,L.Abramov, D.Ornstein and their followers the ergodic theory obtained new perspectives and new links. The idea of entropy as an invariant
of the various objects in the very different part of mathematics became very popular.
Nevertheless, after fifty years of the developing entropy theory there is a room for
generalizations and further application of those ideas. Here we present, the natural generalization of
this notion which is useful even for automorphisms. Moreover, it seems that the entropy as a very special invariant
of automorphism which is very far form others invariants (like spectra) has had up to now no a right framework in the  ergodic theory. We suggested a context - {\it asymptotic invariants of the systems of the compacts with measures} in which entropy and its generalization -{\it scaling entropy} looks in a very natural way. This idea firstly
appeared in the theory of filtrations -decreasing sequences of sigma-fields (or measurable partitions).
We believe that in this terms will be possible to clarify some old problems of the ergodic theory.

\section{Kolmogorov's and Sinai's definitions of entropy. Why
 entropy cannot be deformed?}

We start with discussing an important question related to the notion
of Kolmogorov entropy and its generalizations.
In a descriptive form, the question can be stated as follows: what is an appropriate
framework for considering entropy? and how, within this framework, can one
extend the notion of entropy to the case where the Kolmogorov entropy
of an automorphism vanishes? This problem excited the lively interest of specialists
from the moment A.~N.~Kolmogorov had discovered entropy as an invariant
of dynamical systems. But at that time, this did not lead to any serious consequences.
The attempt to squeeze entropy into the framework of functional analysis,
i.e., into the framework of spectral approach, failed, since the
nature of entropy is obviously non-spectral and even non-operator.
Of course, one can (and did) artificially recast the notion of
entropy in operator terms, but such reformulations do not
give any essentially new information: entropy agrees poorly with traditional
operator considerations, because it is of a completely different nature.
One should search for other appropriate terms and notions.

However, this isolation of entropy is not of course insurmountable. Below
we suggest another framework, important in itself, which naturally embraces both the
Kolmogorov entropy and its generalizations. Our approach is as follows:
given a dynamical system, we consider the {\it associated dynamics
of metrics in a measure space
and its asymptotic invariants}.
Within this approach, it is natural to consider entropy and its generalizations
as {\it one of the basic and simplest asymptotic invariants} of the associated
dynamics of metrics and, consequently, of the dynamical system itself.
But for this we must pass from the Kolmogorov entropy to the
$\varepsilon$-entropy of metric measure spaces and asymptotic characteristics
of their dynamics. This will be done below.

The classical spectral theory of dynamical systems, in particular, the spectral
theory of groups of measure-preserving transformations, studies the associated
dynamics in function spaces, or, from a more modern point of view, in operator
algebras. Our main thesis is as follows:
deep properties of a system are reflected in
how the group acts on the collection of metrics on the phase space, i.e.,
in the associated dynamics of metrics, in particular, in
the asymptotic properties of this dynamics of metrics.
It is this framework that is a natural
place for entropy and its generalizations.

This approach first appeared in connection with the theory of decreasing
sequences of partitions (filtrations)
\cite{V1,V2}. There it turns out to provide a fine classification
of filtrations. The approach is equally fruitful in the theory of dynamical systems
itself.

First recall the basic definitions of entropy theory. The entropy of a discrete
measure
$\mu=(p_1, \dots, p_n)$, $p_i>0$, $\sum p_i=1$, is defined as
$$H(\mu)=-\sum_i p_i \log p_i.
$$

Let $T$ be an automorphism, i.e., a measurable invertible transformation
with an invariant measure, defined on a standard measure  (i.e., in Rokhlin's terms, Lebesgue) space
$(X,{\cal A},\mu)$, where
$\cal A$ is the $\sigma$-algebra of mod~0 classes of all measurable sets
and $\mu$ is a probability measure. Now assume that
$T$ is realized as the right shift in the space of two-sided sequences
of at most countably many symbols. This means that $X$ is the space of
infinite two-sided sequences of symbols
(e.g., $X={\Bbb
N}^{\Bbb Z}$, $T$ shifts a sequence to the right, and
$\mu$ is a shift-invariant (stationary) probability measure). So we
consider a stationary random process. Thus we can obtain (in many ways) a
realization of any automorphism: this is
Rokhlin's theorem on the existence
of a countable generator \cite{R2}. Denote by $\zeta$ the partition
of the space ${\Bbb N}^{\Bbb Z}$ according to the ``past'' of the
process:
an element of $\zeta$ is a class of all sequences with fixed
values of coordinates with negative indices and arbitrary values of coordinates
with nonnegative indices. Let us shift $\zeta$ one place to the right
and consider its average conditional entropy, i.e., the expectation, over all
elements of the past, of the entropy of the conditional distribution of the
zero coordinate given all coordinates with negative indices:
$EH(\zeta \mid T^{-1}\zeta)$.

\medskip\noindent
\textbf{Kolmogorov's theorem.}
{\it The (finite or infinite) nonnegative number
$$EH(\zeta \mid
T^{-1}\zeta)\equiv h(T),$$
called the average conditional entropy (or Shannon information) per step
of a finite-state stationary random process, is an invariant of the automorphism
$T$. In other words, it does not depend on
the particular isomorphic realization of
$T$ as the shift in the space of trajectories of such a process.}
\medskip

Let us return for a while to the history of this discovery.
There are some events related to the statement of this theorem. In his
first paper \cite{K3}, Kolmogorov
interpreted the above statement more widely, but the problem is
that for a continual set of symbols (and even for a countable one, but with an
infinite entropy), the theorem is not true unless we impose some special
conditions on the realization of the automorphism. This was immediately
observed by V.~A.~Rokhlin, who provided a counterexample: an
automoprhism $T$ of algebraic origin and $T$-invariant $\sigma$-algebras
for which the left-hand sides of the above equation are
different for different generated partition  $\zeta$.\footnote{In \cite{K3}, the error was in an illegitimate
passage to the limit along a decreasing sequence of invariant $\sigma$-algebras
(in short, along a decreasing filtration). Curiously enough, the same error was
made by N.~Wiener in an important passage of his well-known book on
the nonlinear theory of random processes, as well as by many other authors.
The point is that along increasing sequences of $\sigma$-algebras (in short,
along increasing filtrations), the passage to the limit is obviously legitimate,
and this provokes one to assume that the same is true for decreasing filtrations.
However, the theory of decreasing filtrations, and especially their classification,
is much finer and more interesting than the theory of increasing filtrations.
As a rule, the passage to the limit along a decreasing filtration
is not possible (for the theory of decreasing filtrations, see
\cite{V2}).} In his second note
\cite{K4}, A.~N.~Kolmogorov corrected the statement by imposing {\it a priori}
conditions on automorphisms in the definition of entropy.
However, for Bernoulli schemes, i.e., sequences of independent random variables
with a finite or countable state space, as well as for many other cases,
the invariant was well-defined already in the first paper
\cite{K3}. For the above statement to be true in full generality, one needed
theorems on generators and special invariant $\sigma$-algebras which were not
yet known. The above-mentioned theorem on the existence of a countable
generator for any aperiodic automorphism was proved  somewhat later by V.~A.~Rokhlin:
any automorphism can be realized as the shift in the space of trajectories of a process
with at most countably many states. This recovered the generality
of Kolmogorov's theorem. Somewhat later, W.~Krieger proved the existence of
a finite generator for automorphisms with finite entropy. The same fact, but with
a worse estimate, was proved by A.~N.~Livshits
(1950--2008).

A much simpler definition, generally accepted nowadays, was suggested
by Ya.~G.~Sinai soon after the appearance of Kolmogorov's work. It is not related to information,
but is rather of geometric, or
combinatorial, nature.

\medskip\noindent
\textbf{Sinai's theorem {\rm\cite{S1}}.}
{\it Let $T$ be an automorphism of a standard measure space
$(X,{\cal A},\mu)$ and $\xi$ be a finite measurable partition of this space.
Let $T^0\xi=\xi, T\xi, T^2\xi, \dots, T^{n-1}\xi$ be the successive
$T$-images of $\xi$ and
$\xi^n_T=\bigvee_{i=0}^{n-1} T^i\xi$ be the product of the first $n$ images.
Then the following (finite or infinite) limit does exist:
$$\lim_{n\to \infty} \frac{H(\xi_T^n)}{n}\equiv h(T,\xi);$$
the expression
$$\sup_{\xi}h(T,\xi)=h(T),$$
where the supremum is taken over all finite partitions of $X$,
coincides with the entropy of $T$ defined above
(or may be taken as its definition).}
\medskip

Kolmogorov's and Sinai's definitions are of different nature and
have different interpretations (see below); their equivalence is not
quite obvious. An important fact: the expression
$h(T,\xi)$ is continuous in
$\xi$ on the space of all finite partitions endowed with the so-called
entropy metric. It is this continuity that allows one to effectively compute
the entropy via approximations. As mentioned above, the positivity of entropy
distinguishes an important class of automorphisms whose properties differ
from those of zero-entropy automorphisms.

Assume that the entropy of an automorphism is zero. Then the sequence
$H(\xi_T^n)$ grows sublinear. Can we replace
the linear scaling by another one so that to obtain a new invariant of the automorphism? In other
words, can we deform Sinai's definition in such a way that the new invariant would
distinguish at least some zero-entropy automorphisms?
It turns out that Kolmogorov's entropy in Sinai's definition cannot be
deformed in the following literal sense.

\begin{theorem} For every ergodic automorphism $T$ and any sequence of
positive numbers $ \{c_n\}$ satisfying the condition
$\lim_n \frac{c_n}{n}=0$, there exists a finite partition $\xi$ such that
$$\lim_n \frac{H(\xi^n_T)}{c_n}=+\infty.$$
\end{theorem}

In other words, no striclty sublinear growth of the entropy
$H(\xi_T^n)$ provides a new invariant. This effect is due to the fact that
the partition $\xi^n_T$ in Theorem~1 can be very fine on a set of small measure,
thus leading to an artificially high value of the entropy. Not quite accurately,
one can formulate this as follows: on a set of small measure, the growth
of the entropy can be almost linear (note that it cannot be superlinear).
An analogous observation explains why a similar analog of the original Kolmogorov's
definition does not lead to reasonable new invariants.
This is a manifestation of the general principle of ergodic theory expressed
by Rokhlin's lemma: ergodic automorphisms are indistinguishable up to
a set of small measure.

Nevertheless, one can refine the idea of scaling the growth. For this, one
should consider entropy only up to small changes in the partition
$\xi_T^n$, in other words, involve the notion of $\varepsilon$-entropy.
First we will describe new invariants in terms as close as possible
to the traditional ones, i.e., using partitions as in Sinai's definition, and
then turn to the richer language of metrics.

\section{The $\varepsilon$-entropy of a measure space and the definition
of scaling entropy}

Consider the following function of a partition $\xi$ and a positive number
$\varepsilon$:
$$H_{\varepsilon}(\xi)=\inf_{A:\;\mu A>1-\varepsilon} H(\xi\mid_A).$$
By $\xi\vert_A $ we mean the partition of a set
$A$ of positive measure, with the restriction
of the measure $\mu$ to $A$ renormalized to unity, whose elements
are the intersections of the elements of
$\xi$ with $A$. Observe that this function monotonically
decreases as $\varepsilon$ grows, taking the value
$H(\xi)$ for
$\varepsilon=0$ and vanishing for  $\varepsilon=1$. We will use it to
define
the scaling entropy of an automorphism. Consider the function
$$H_{\varepsilon}(\xi_T^{n}),$$
which depends on $n$, $\varepsilon >0$, and
a partition $\xi$, and let us study its growth.

\begin{definition}[\cite{V4}]
We say that a sequence of positive numbers
$\{c_n\}$ is a scaling sequence for an ergodic transformation $T$ if for every
finite partition $\xi$
$$\lim_{\varepsilon \to 0} \limsup_n \frac{H_{\varepsilon}(\xi_T^{n})}{c_n}<\infty,$$
and there exists a partition $\xi$ such that
$$\lim_{\varepsilon \to 0} \liminf_n \frac{h_T(\xi n, \varepsilon)}{c_n} > 0. $$
(Note that all such sequences $\{c_n\}$ are equivalent as $n\to \infty$.)
\end{definition}

\begin{theorem}
The class of scaling sequences for a given transformation is a metric invariant of
the transformation. This invariant distinguishes transformations with
zero Kolmogorov entropy. The class of sequences
$\{c_n\} \sim \{n\}$ corresponds to the Kolmogorov entropy; with this
scaling, the function
$H_{\varepsilon,n}(\xi)$ for small $\varepsilon$ does not depend on $\varepsilon$.
\end{theorem}

Sometimes, in a class of equivalent sequences one can choose a sequence
suitable for all transformations having this class as the scaling one.
Then the invariant we obtain is not merely a class, but a number, called
the {\it scaling entropy}. It is not clear whether one can always make such
a choice. In \cite{KT}, the so-called ``slow'' entropy was introduced
to distinguish actions of the group
${\Bbb Z}^2$; it resembles our
scaling entropy.

Let us give several examples.

\smallskip\noindent{\bf Example 1.}
If the scaling sequence is linear, i.e.,
$\{c_n\} \sim \{n\}$, then we obtain the Kolmogorov entropy. Of course,
this growth is the maximum possible one for the group
$\Bbb Z$. In the case of a zero Kolmogorov entropy, the scaling is sublinear.

\smallskip\noindent{\bf Example 2.} The scaling  class corresponding to transformations with discrete spectrum
is that of bounded sequences:
$\sup c_n <\infty$. Thus for measure-preserving isometries, i.e., shifts on
compact groups, the scaling entropy vanishes for any scaling.
In a slightly different formulation, this fact was first observed by
S.~Ferenczi \cite{Fe}. In appearance, it resembles A.~Kushnirenko's
\cite{Ku} result on the Kirillov--Kushnirenko entropy, also called
{\it sequential entropy}: the class of automorphisms for which the sequential
entropy vanishes for any sequence coincides with the class of automorphisms with
discrete spectrum.

However, this resemblance is superficial, because the notion of scaling
entropy differs crucially from that
of Kirillov--Kushnirenko entropy.

\smallskip\noindent{\bf Example 3.}
In exactly the same way as above, a scaling sequence can be also defined
for flows. In two independent articles,
A.~Kushnirenko's paper \cite{Ku} mentioned above  and M.~Ratner's paper \cite{R}
on horocycle flows on surfaces of negative curvature, it was proved that
different Cartesian powers of such flows are nonisomorphic.
The distinguishing invariant in
\cite{Ku} was the sequentional entropy, and in
\cite{R} M.~Ratner used an idea of J.~Feldman \cite{F} and constructed an
invariant that resembles a special example of scaling entropy. Comparing the
invariant from \cite{R} with our definition, one can conjecture that the scaling
sequence for the $k$th power of a horocycle flow is logarithmic:
$\{c_n\} \sim \{(\log n)^k\}$.

\smallskip\noindent{\bf Example 4.}
A challenging problem is to find scaling sequences for
adic automorphisms, e.g., the Pascal or
Young automorphisms. Presumably, these automorphisms have singular continuous
spectra.  If the scaling sequence is
not bounded,
it would follow from above (see Example~2) that
the spectrum is not purely discrete; this problem is still open.

\section{Admissible metrics instead of partitions}

The approach to scaling entropy described above needs further development.
Instead of using the theory of partitions, which is a traditional
approach in ergodic theory,
one should involve the more flexible techniques of metrics and metric spaces,
which turn to be useful in many problems of measure theory. Below we illustrate this with the
example of the theory of scaling entropy, which includes the ordinary entropy
theory. In the author's opinion, this approach must also be fruitful in
applications to the general isomorphism problem in ergodic theory.

Any finite or countable measurable partition $\xi$ determines a semi-metric:
$$\rho(x,y)= \delta_{C(x),C(y)},$$
where $C(x)$ stands for the element of $\xi$ that contains $x$.

Thus we can replace manipulations with measurable partitions by the analysis
of the corresponding semi-metrics; in other words, the transition to
(semi)metrics tautologically includes the theory of partitions as a special case.
But considering general metrics and semi-metrics also opens up new possibilities.

Our approach is as follows:
{\it instead of studying the collection of Borel measures on a fixed
metric (or topological) space, which is a usual practice, we consider a set of metrics on a fixed measure
space}. From the viewpoint of ergodic theory and probabilistic considerations, such a shift
is quite natural. Let us introduce the notion of an
{\it admissible (semi)metric} on a Lebesgue space
$(X,\mu)$.

\begin{definition}
A (semi)metric $\rho$ on a Lebesgue space $(X,\mu)$ with continuous measure
is called admissible if the following conditions are satisfied:

{\rm1)} The function $\rho(\cdot,\cdot)$  is a measurable
function of two variables, defined on
a set of full measure (depending on the metric) in the Cartesian square
of the space $(X,\mu)$, that satisfies the axioms of a (semi)metric on this set.

{\rm2)} The space $(X,\rho)$, regarded as a (semi)metric space, is quasi-compact,
i.e., it turns into a compact space after taking the quotient by the equivalence relation
$x\sim y \Leftrightarrow \rho(x,y)=0$.
\end{definition}

Thus a well-defined notion is not an individual (semi)metric on a Lebesgue
space, but a class of  mod~0 coinciding (semi)metrics, so that one should
speak about classes of mod~0 coinciding (semi)metric spaces.
As a rule, checking that assertions under consideration are well-defined with respect to mod~0
equivalence presents no problem. Nevertheless, there are some
subtleties, e.g.,
in the understanding of the triangle inequality (it should hold
for all triples of points from the set of full measure on which the metric
is defined, and not for almost all triples of points). The admissible
(classes of) (semi)metrics form a convex cone
$\cal R$ in the space of measurable nonnegative functions of two variables
on the space $(X,\mu)$. We call $\cal R$ the cone of (classes of) admissible
(semi)metrics; it is closed under supremum of finitely many metrics.
This is a canonical object, provided that we restrict ourselves to Lebesgue spaces
with continuous measure. The geometry of this cone is of great interest;
it is poorly studied.

\section{The $\varepsilon$-entropy of measures in metric spaces}

The following definition of the $\varepsilon$-entropy of measures in metric spaces
is also essentially due to A.~N.~Kolmogorov (see \cite{K1}). We change only
one detail, which is not very important; namely, we estimate the closeness of measures
in the Kantorovich metric rather than by counting the number of points
in an $\varepsilon$-net.

\begin{definition}
Let $\mu$ be a Borel probability measure on a separable metric space
$(X,\rho)$. Define a function $H(\rho, \mu, \varepsilon)$
as follows: $$H(\rho, \mu, \varepsilon)=\inf \{H(\nu):
k_{\rho}(\mu ,\nu)<\varepsilon\},$$
where $\nu$ ranges over the set of discrete measures and
$k_{\rho}$ is the Kantorovich distance between measures in the metric space
$(X,\rho)$.
\end{definition}

Recall the definition of the Kantorovich (transportation) metric
\cite{Kan} on the space of measures defined on a compact metric space.
\footnote{The necessity to use Kantorovich metric on the set of the probability measures on the metric space,has two explanations: the first --- this metric is an analog $L^1$-metric for measures which is important in what follows, and  secondly, Kantorovich metric is maximal among all metrics on the set of probability measures on the metric space with property $d(\delta_x,\delta_y)=\rho(x,y)$ (see\cite{M}).}

Let $(X,\rho)$ be a compact metric space, and let
$\mu_1,\mu_2$ be two Borel probability measures on
$X$. Then
$$
k_{\rho}(\mu_1,\mu_2)=\inf_{\Psi}\int\int_{X\times X} \rho(x,y)\,d\Psi(x,y),
$$
where the infimum is taken over all probability measures
$\Psi$ on the space $X\times X$ whose projection to the first coordinate
coincides with $\mu_1$ and projection to the second coordinate coincides with
$\mu_2$. In other words, $\Psi$ ranges over the set of measures with
given marginal projections.

Note that the above definition makes sense also in the case where
the space is not compact, because a probability measure in a complete separable
space is supported, up to any positive
$\varepsilon$, by compact subsets. In the case of semi-metrics, the definitions
also remain meaningful. If we are given a finite partition
$\xi$ of a measure space $(X,\mu)$, then its entropy coincides with
the $\varepsilon$-entropy (for sufficiently small
$\varepsilon$) of the space $(X/\xi, \mu/\xi)$
of elements of $\xi$ with the discrete (semi)metric
$\rho_{\xi}$: $$ H(\rho_{\xi}, \mu, \varepsilon)= H_{\varepsilon}(\xi).$$

\section{Dynamics of metrics in a measure space as an appropriate
framework for entropy}

The classical functional analysis suggests to consider, instead of various
objects, the spaces of functions on these objects. The spectral theory of
dynamical systems is the result of following this recommendation: instead of
a transformation of the phase space one considers a unitary operator
in the corresponding $L^2$ space. But somehow
these considerations  have been hitherto limited
to functions of one variable running over the phase space of the
system. However, one may consider actions of Cartesian powers of the dynamical system
in spaces of functions of several variables, while preserving the
separation of variables;
for instance,
in the space of functions of two variables, namely, on the cone
of admissible metrics. Clearly, in this way we obtain much more information
on the system than when considering an action in the space of functions of one
variable, and thus increase the possibilities of analyzing the properties
of dynamical systems.\footnote{Of course, Cartesian powers are widely used
in ergodic theory, but usually one considers the Cartesian square merely as
an automorphism of a measure space, not fixing the structure of a direct product.}
This leads us to new interesting and important problems.

Let  $\rho$  be a (semi)metric and $T$ be an automorphism. Denote by
$\rho_T$ the (semi)metric $\rho_T(x,y)=\rho(Tx,Ty)$. The image of an admissible
metric is an admissible metric. Thus there is a natural action of
the group of measure-preserving transformations on the cone  $\cal R$.

Our main thesis is as follows: it is the
{\it asymptotic theory of iterations of metrics in a space with a fixed measure}
under automorphisms that is an appropriate framework for considering
both Kolmogorov and scaling entropies (and their generalizations),
as well as other invariants of automorphisms.

Given an admissible metric $\rho$ on a space $(X,\mu)$
and an automorphism $T$, we can construct a sequence of new metrics. For this,
we must take the orbit of $\rho$ in the cone of admissible metrics under the
action of $T$ and then form symmetric combinations of the first
$n$ elements of the orbit. The following two sequences of metrics associated with
a given metric $\rho$ and a given automorphism $T$ are especially important:
the {\it uniform metric}
$$
\rho^n_T=\sup_{i=0,\dots n-1} \rho_{T^i} \quad\text{where}\quad \rho_{T^i}(x,y)=\rho(T^ix,T^iy)
$$
(in the  ordinary setting, it corresponds to the product of partitions:
$\rho_{\xi^n_T}=\sup_{i=0,\dots n-1}\rho_{T^i\xi}$)
and the
{\it average metric}
$$\widehat{\rho}^{\,n}_T=\frac{1}{n}\sum_{i=0}^{n-1}\rho_{T^i}.$$

Below we restrict ourselves to the first of them, nevertheless we  believe that average
metric must be the main. It has no interpretation in terms of partitions because
average of partition has no sense. But from the technical point of view it is much more convenient.
because it has {\it ergodic feature} \footnote{It make sense to consider quadratic and other average:
 $$[\widehat{\rho}^{\,n}_T]^p=\frac{1}{n}\sum_{i=0}^{n-1}\rho_{T^i}^p,$$
 $p=2$ is useful in the case Riemannian space.}

{\it The problem we set out is to study the asymptotic behavior of this sequence of metrics
as $n$ goes to infinity.}

\medskip\noindent
\textbf{A scaling sequence.} Consider the growth of the $\varepsilon$-entropy
of a measure $\mu$ in the sequence of compact metric spaces
$(X, \rho^{\,n}_T)$ and introduce the class of monotone sequences
of positive numbers $\{c_n\}$ such that
$$0<\lim_{\varepsilon\to 0}\liminf_n\frac{H(\rho^n_T,\mu,\varepsilon)}{c_n}\leq
\lim_{\varepsilon\to 0}\limsup_n\frac{H(\rho^n_T,\mu,\varepsilon)}{c_n}<\infty.$$
We say that a sequence ${c_n}$ from this class is a {\it scaling sequence
for the automorphism $T$ and the metric $\rho$}. This normalization corresponds to the one
we considered above when defining the scaling entropy of an automorphism
via partitions.

\begin{theorem}[On the scaling entropy of an ergodic automorphism]
For an ergodic automorphism $T$, the class of scaling sequences $\{c_n\}$ does not depend on the choice
of a metric $\rho$ in the class of metrics satisfying, along with
conditions {\rm1)} and {\rm2)} above, the following condition:

{\rm3)}The metric (not semi-metric!) is a generic metric in the sense of topology in the
space of admissible metric.
\end{theorem}

If for some canonical choice of a sequence from the equivalence class,
the limit
$$\lim_{\varepsilon\to0}\lim_n\frac{H(\rho^n_T,\mu,\varepsilon)}{c_n}$$
does exist, we call it  the
{\it scaling entropy} of $T$, indicating the scaling sequence.

\medskip
\textbf{Problem.}  Does the class of scaling sequences change when we substitute
in its definition the metrics $\rho^{\,n}_T$ with the average metrics $\widehat{\rho}^{\,n}_T$?
\medskip

Thus the scaling entropy (and, in particular, the Kolmogorov entropy) is a natural
asymptotic invariant of a sequence of compact measure spaces. The methodological
advantage of passing from partitions to continuous (semi)metrics is that
passing to the limit and taking the supremum over all finite partitions
is now replaced by considering an appropriate semi-metric.
But it is much more important that now the original problem reduces to a circle
of asymptotic questions concerning the behavior of a sequence of
metrics. And the scaling
entropy is merely one of the asymptotic characteristics of
these metrics (the coarsest one). It characterizes the growth of the ``dimension'' of the compact space,
or, in other words, the {\it asymptotics of its Hausdorff dimension}.
It is the asymptotics of the dynamics of the metric measure spaces
$(X, \mu, \rho^n_T)$
that is the appropriate framework for considering entropy and its generalizations
we mentioned above.

Let us briefly describe the program we propose. Consider any admissible metric
on a Lebesgue space $(X,\mu)$ with a fixed continuous measure and a
measure-preserving automorphism $T$ (or a group of automorphisms $G$).
We suggest to study asymptotic invariants of the sequence of metrics
$\rho^n_T$ introduced above. The asymptotic characteristics of this sequence
do not depend on the original metric and thus characterize only invariant
properties of the automorphism.

As $n$ grows, the metrics change in a quite complicated way, but presumably
there is a number of coarse asymptotic invariants such as entropy. The scaling entropy
is a simplest asymptotic invariant, which describes the growth of the cardinality
of $\varepsilon$-nets, or the growth of the Hausdorff dimension of the compact
space. More involved asymptotic invariants characterize not only the asymptotics
of individual compact spaces, but also the asymptotics of their relative position.
It is not yet known whether for the sequence of
compact spaces  there exists a limit object
(in the example with filtrations considered below, such a limit object
does exist). It may happen that in this setting
limit objects also do exist and can be characterized
in a more or less explicit way. Presumably, the study of such objects would
allow one to solve the isomorphism problem for automorphisms with completely
positive entropy.

\section{Relation to invariants of metric triples and their dynamics}

Let us relate our considerations to the theory of metric triples, or Gromov triples,
or $mm$-spaces in Gromov's terminology. Recall that M.~Gromov
\cite{G}\footnote{We quote from the previous edition of the book,
published in 2001.} suggested a complete invariant of triples
$(X,\rho,\mu)$ with respect to $\mu$-preserving isometries of the space $(X,\rho)$.
Here $X$ is a space, $\rho$ is a metric on $X$
that turns it into a Polish space, and $\mu$ is a
fully supported Borel probability measure.
In the formulation due to the author of the present paper
(see \cite{V3}), this invariant looks as a probability measure on the cone of
nonnegative infinite symmetric matrices with countably many rows and columns
satisfying the triangle inequality, i.e., on the cone of so-called distance
matrices. This invariant can be interpreted as a random (semi)metric on the
positive integers, or as a random distance matrix. A large variety of known
invariants of metric triples can easily be computed via this invariant, i.e,
via a random matrix. For example, in these terms one can easily describe the
$\varepsilon$-entropy of the triple (or, say it another way, of the measure
$\mu$ in the metric space $(X,\rho)$).

In the above scheme, we considered
a sequence of metrics on the same measure space, i.e.,
a sequence of metric triples that differed only
by metrics. And the problem was to find asymptotic invariants of these triples.
The complete invariant of triples we have just described
allows one to easily compute
the $\varepsilon$-entropy, and thus to find the asymptotics
of the $\varepsilon$-entropies of the sequences of triples, which, in our
definition, is exactly the scaled entropy. Note that, as mentioned above,
the result does not in fact depend on the original metric. One can also suggest
other coarse asymptotic characteristics of sequences of triples; however,
the choice of an appropriate characteristic should be determined by the problem under consideration.

For example, how  can one formulate Ornstein's results on the classification and
characterization of Bernoulli automorphisms
(via the  $\bar d$-metric or the $VWB$-property)
as an assertion about some special asymptotic type of a sequence of
metric triples? According to Ornstein, each such asymptotic type is determined
by just one positive number, the entropy of the Bernoulli automorphism. In
other words, the problem is to find a geometric description of the Bernoulli type
of sequences of metrics. However, this problem is apparently still far from being solved.
One may hope that examples of non-Bernoulli automorphisms with completely positive
entropy would be more fully explained in terms of the asymptotic dynamics
of metric spaces.

\section{A parallel with the theory of filtrations and the dynamics of iterations
of the Kantorovich metric}

Another dynamics of metrics, which had appeared much earlier, is related to the theory
of filtrations and, in particular, to the entropy of filtrations. Though this
dynamics is more complicated, an attractive feature of this approach is the existence of limit
objects; this allows one to develop the study much further than is currently
done for the project described above. Let us briefly mention some definitions
and examples.

A filtration is a decreasing sequence of
$\sigma$-algebras. An example of a filtration is the sequence of pasts of a random process
$\{T^{-n}{\cal A}\},$ $n=0,1,\dots$, where $\cal A$
is a $T$-invariant $\sigma$-algebra (i.e., $T^{-1}{\cal A}\subset \cal A$).
A filtration is called ergodic if the intersection
of the $\sigma$-algebras $T^{-n}{\cal A}$ is the trivial algebra.

It $T$ is a one-sided Bernoulli shift, then the corresponding filtration of pasts
(which is ergodic by Kolmogorov's zero-one law) is called {\it standard}.
There exist ergodic filtrations whose finite parts are isomorphic to Bernoulli
filtrations but that are not isomorphic to Bernoulli filtrations as a whole.
The following {\it dynamics of (semi)metrics generated by a filtration}
is of fundamental importance.
It is more convenient  to pass from a filtration of
$\sigma$-algebras to a (decreasing) filtration of partitions
$\xi_1 \succ \xi_2 \succ \dots\,$.
Suppose we have a filtration $\xi_n$ with trivial intersection - $\bigcap_n \xi_n=\nu$
($\nu$ is trivial partition). Consider an admissible metric $\rho$ on a measure space
$(X,\mu)$ and construct the sequence of metrics $\rho_n$ in the different way than we had for the
case of automorphisms in the previous paragraph. Namely we use the Kantorovich iterations
constructed with the help of the filtration.

Let
\begin{eqnarray*}
\rho_0=\rho,\quad \rho_1(x,y)=k_{\rho_1}(\mu^{C_1(x)},\mu^{C_1(y)}), \dots, \\
\rho_m(x,y)=k_{\rho_{m-1}}(\mu^{C_m}(x),\mu^{C_m}(y)),\dots,
\end{eqnarray*}
where $C_m(x)$ is the element of $\xi_m$ that contains $x$,
$\mu^C$ is the conditional measure on an element $C$,
and $k_{\rho}$ is the Kantorovich distance between measures on the metric measure
space $(X,\rho)$. In other words, the distance between points
$x$ and $y$ in the $n$th semi-metric is the Kantorovich distance between the
conditional measures on the elements
$C_n(x)$ and $C_n(y)$ of the $n$th partition
with respect
to the $(n-1)$th semi-metric.

We obtain a sequence
$\{\rho_m\}_0^{\infty}$ of semi-metrics on the space
$(X,\mu)$. In terms of the asymptotic behavior of this sequence,
one can express many (possibly all) invariant asymptotic properties of the
filtration. Remarkably, these asymptotic properties do not depend on the
choice of the original metric from a very wide class of metrics --- exactly as in
the program considered above.

The main example is related to the standardness criterion, see the author's papers
\cite{V1,V2}.

\begin{theorem}
The sequence of iterated metrics tends to a degenerate metric (i.e., the metric
space contracts to a point) if and only if the filtration is standard (Bernoulli).
\end{theorem}

This means that the scaling entropy of the sequence of metrics vanishes
for any increasing sequence
$\{c_n\}$. One may compare this fact with the dynamics of metrics under
automorphisms with discrete  spectrum (see above).

\bigskip
\textbf{EXAMPLE: scaling entropy for the random walk in a random environment.}

Let us give a newer example. Consider a random walk in a random environment: namely,
the simple random walk (and the corresponding Markov process)
on the set of all $\{0,1\}$-configurations of the lattice ${\Bbb Z}^d$
equipped with Bernoulli measure $(1/2,1/2)$.

For $d=1$, it is so called $(T,T^{-1})$ transformation, where $T$ is Bernoulli automorphism;
the Markov shift in this case is a non-Bernoulli $K$-automorphism even not loosely Bernoulli -
(S.~Kalikow \cite{Ka}). As F.~Hollander and J.~Steif \cite{STHO} shown the same holds for
$d=2$;  for $d>2$ the Markov shift is already Bernoulli, but in the case $d=3,4$ the natural generator
is very weak Bernoulli but not weak, and is a weak Bernoulli for $d>4$.
All this results simulated the question: what can be said about the filtration
of pasts of these processes? The conjecture that in the case $d=1$ the filtration of the past is not standard
 was formulated by A.Vershik in 70-th. Now we can answer on properties of the filtration of the pasts
 and calculate scaling entropy.
\begin{theorem}
The scaling entropy for the filtration of the past in the case of
random walk in the space of configuration on the lattice ${\Bbb Z}^d$ of dimension $d$
is normalized by sequences which are equivalent to $c_n=n^{\frac{d}{2}}$.
\footnote{The growth $c_n=n^k$ in this situation corresponds to a logarithmic growth for
the group $\Bbb Z$ $(\ln)^k$, since here we consider the group $\sum
{\Bbb Z}_{2d}$, which has infinitely
many generators, and has exponential growth of the number of words, so usual (Kolmogorov) entropy
for this group must be normalized with the growth $c_n=(2d)^n$; the result above
can be compared with growth $c_n=(\ln n)^{d/2}$
 for the group $\Bbb Z$, which is the same as in the example of horocycle flow in the examples 3 above.}
\end{theorem}
\medskip
 D.~Heicklen and C.~Hoffman \cite{HH} proved
the nonstandardness for the dimension
$d=1$, and then essentially computed the entropy for
$d=1$. A.~Vershik and A.~Gorbulsky in
\cite{VG} proved
the nonstandardness for $d>1$, and, found the scaling entropy in general case.
It showed, by the way, that the filtrations are nonisomorphic for different $d$. The latter result implies
that if the dimensions $d$ of lattices are different, than the Markov
process of the random walk on one lattice cannot be encoded in an invertible way
into the shift on the other lattice, though (for $d>3$)
all these Markov shifts are Bernoulli.

The dynamics of metrics in the case of filtrations is closely related to the
construction of the so-called tower of measures, which allows one to construct
limit metric spaces for a sequence of iterated compact spaces
(see \cite{V2,V4,VG}). The corresponding combinatorics is quite interesting and
has relevance to actions of groups of automorphisms of trees and close groups.

\bigskip
Translated by N.~V.~Tsilevich.


\begin{thebibliography} {99}
\bibitem{F}
J.~Feldman, $r$-entropy, equipartition, and Ornstein's isomorphism
theorem in ${\bf R}^n$. {\it Israel Math. J.} \textbf{36} (1980),
321--345.

\bibitem{Fe}
S.~Ferenczi, Measure-theoretic complexity of ergodic systems. {\it Israel Math. J.}
\textbf{100} (1997), 180--207.

\bibitem{G}
M.~Gromov, {\it Metric Structures for Riemannian and Non-Riemannian
Spaces}. Birkh\"auser Boston,
Boston, MA, 2007.

\bibitem{HH} D.~Heicklen and  C.~Hoffman,  $T,T^{-1}$ is not
standard. {\it Ergodic Theory Dynam. Systems} \textbf{18} (1998), No.~4,
875--878.

\bibitem{STHO} F.~Hollander and J.~Steif, Random walk in random
scenery. {\it IMS Lect. Notes} \textbf{48} (2006), 53--65.

\bibitem{Ka}
S.~A.~Kalikow, $T,\,T\sp{-1}$ transformation is not loosely
Bernoulli.  {\it Ann. of Math. (2)}  \textbf{115} (1982), No.~2, 393--409.

\bibitem{Kan} L.~V.~Kantorovich, On the translocation of masses. {\it
Dokl. Akad. Nauk SSSR}
{\bf 37} (1942), No.~7--8, 227--229.

\bibitem{KT}
A.~Katok and J.-P.~Thouvenot, Slow entropy type invariants and
smooth realization of commuting measure preserving transformation.
{\it Ann. Inst. H. Poincare} \textbf{33} (1997), 323--338.


\bibitem{K1} A.~N.~Kolmogorov, Information transmission theory. In: {\it Session
of the Academy of Sciences of the USSR on Scientific Problems of Industrial Automation},
October 15--20, 1956. Izdat. Akad. Nauk, 1957,
pp.~66--99.

\bibitem{K3} A.~N.~Kolmogorov, A new metric invariant of
transient dynamical systems and automorphisms in Lebesgue spaces. {\it
Dokl. Akad. Nauk SSSR} {\bf 119} (1958), 861--864.

\bibitem{K4} A.~N.~Kolmogorov,
Entropy per unit time as a metric invariant of automorphisms. {\it
Dokl. Akad. Nauk SSSR}  {\bf 124}  (1959), 754--755.


\bibitem{Ku} A.~G.~Kushnirenko,
Metric invariants of entropy type. {\it Uspekhi Mat. Nauk} {\bf 22}  (1967),
 No.~5 (137), 57--65.

\bibitem{M} J.~Melleray, F.~Petrov, A.Vershik, Linearly rigid metric space
and the embedding problem. {\it Fund.Math.} {\bf 199} (2008), No.2, 177-194.

\bibitem{R}
M.~Ratner, Some invariants of Kakutani equivalence. {\it Israel Math.
J.} \textbf{38} (1981),  232--240.

\bibitem{R2} V.~A.~Rokhlin,
Lectures on the entropy theory of transformations with invariant measure. {\it
Uspekhi Mat. Nauk}  {\bf 22}  (1967),  No. 5 (137), 3--56.

\bibitem{S1} Ya.~G.~Sinai,
On the concept of entropy for a dynamic system. {\it Dokl. Akad. Nauk SSSR}
{\bf  124}  (1959), 768--771.

\bibitem{S2} Ya.~G.~Sinai,
About A.~N.~Kolmogorov's work on the entropy of
dynamical systems. {\it Ergodic Theory Dynam. Systems} {\bf 8} (1988),
501--502.

\bibitem{V1} A.~M.~Vershik, Decreasing sequences of measurable partitions
and their applications, {\it Dokl. Akad. Nauk SSSR},
 {\bf 193} (1970), No.~4, 748--751.

\bibitem{V2} A.~M.~Vershik,
Theory of decreasing sequences of measurable partitions.
{\it St.~Petersburg Math. J.} {\bf 6} (1995), No.~4, 705--761.

\bibitem{V4} A.~M.~Vershik,
Dynamic theory of growth in groups: entropy, boundaries, examples. {\it
Russian Math. Surveys} {\bf 55}  (2000),  No. 4, 667--733.

\bibitem{V3} A.~M.~Vershik,
Random metric spaces and universality. {\it  Russian Math. Surveys} {\bf 59}  (2004),  No.~2, 259--295.

\bibitem{VG}   A.~M.~Vershik and A.~D.~Gorbulsky,
The scaling entropy of filtrations of $\sigma$-algebras, {\it Theory of prob and appl.}
{\bf 52} (2007) No.3, 446-467.

\bibitem{E} {\it Mathematical events of the XX century: look from St.Petersburg}.
Ed. A.Vershik. MCCME, Moscow 2009.(In Russian),

\end{thebibliography}
\end{document}